\documentclass[11pt, letterpaper, oneside]{amsart}

\usepackage{latexsym, amsmath, amssymb, amsfonts, amscd}
\usepackage{amsthm}
\usepackage{t1enc}
\usepackage[mathscr]{eucal}
\usepackage{indentfirst}
\usepackage{graphicx, pb-diagram}
\usepackage{fancyhdr}
\usepackage{fancybox}
\usepackage{enumerate}
\usepackage[all]{xy}
\usepackage{url}

\theoremstyle{plain}
\newtheorem{thm}{Theorem}[section]
\newtheorem*{thmno}{Theorem} 

\newtheorem{prop}[thm]{Proposition}
\newtheorem{lemma}[thm]{Lemma}
\newtheorem{cor}[thm]{Corollary}

\theoremstyle{definition}
\newtheorem{defi}[thm]{Definition}

\theoremstyle{remark}
\newtheorem{remark}[thm]{Remark}
\newtheorem{ep}[thm]{Example}

\newcommand{\ZZ}{\ensuremath{\mathbb Z}}
\newcommand{\CC}{\ensuremath{\mathbb C}}
\newcommand{\RR}{\ensuremath{\mathbb R}}
\newcommand{\kk}{\ensuremath{\frak{k}}}
\newcommand{\g}{\ensuremath{\frak{g}}}
\newcommand{\h}{\ensuremath{\frak{h}}}

\newcommand{\sll}{\ensuremath{\frak{sl}}}
\newcommand{\su}{\ensuremath{\frak{su}}}
\newcommand{\so}{\ensuremath{\frak{so}}}
\newcommand{\spp}{\ensuremath{\frak{sp}}}

\newcommand{\cL}{\mathcal{L}}
\newcommand{\cO}{\mathcal{O}}

\newcommand{\sh}{\sharp}

\begin{document}

\title{A construction for coisotropic subalgebras of Lie bialgebras}
\author{Marco Zambon}
\address{Departamento de Matem\'aticas,
Universidad Aut\'onoma de Madrid,
Campus de Cantoblanco,
28049 - Madrid, Spain}
\email{marco.zambon@uam.es}

\date{}

\begin{abstract}
Given a Lie bialgebra $(\g,\g^*)$,
we present an explicit procedure to construct coisotropic subalgebras, i.e.  Lie subalgebras of $\g$ whose annihilator is a Lie subalgebra of $\g^*$.  
We write down families of examples
for the case that $\g$ is a classical complex simple Lie algebra.
\end{abstract}

\maketitle
\tableofcontents

\section{Introduction}\label{intro}
A \emph{Lie bialgebra} \cite{Dr6871} structure on a Lie algebra $(\g ,[\bullet,\bullet])$ is  a  degree 1 derivation $\delta$ of $\wedge^{\bullet}\g$ which squares to zero and 
satisfies $\delta([X,Y])=[\delta(X),Y]+[X,\delta(Y)]$.
Dualizing $\delta|_{\g}\colon \g\rightarrow \wedge^2 \g$ one obtains 
a Lie bracket on $\g^*$, encoding $\delta$, so that the Lie algebra structures on $\g$ and $\g^*$ are compatible.
The aim of this paper is to construct Lie subalgebras $\h$ of $\g$ with the property that $\h^{\circ}$, the subspace of $\g^*$ consisting of elements that vanish on $\h$, is a Lie subalgebra of $\g^*$. Such an $\h$ is called \emph{coisotropic subalgebra}.

Our main result (Thm. \ref{rmat}) is a explicit and computationally friendly construction that works for Lie bialgebras arising from  $r$-matrices. Recall that any $r$-matrix on a Lie algebra $\g$, i.e. any $\pi\in \wedge^2\g$ such that $[\pi,\pi]$ is $ad$-invariant, 
gives rise to a Lie bialgebra by setting $\delta=[\pi,\bullet]$.
Our result can be phrased as follows:
\begin{thmno}
Let $\g$ be a Lie bialgebra arising from an  $r$-matrix $\pi$.
Suppose $X\in \g$ satisfies
\begin{equation*} 
[X,[X,\pi]]=\lambda[X,\pi] \text{ for some }\lambda\in \RR. \end{equation*} Then
 the image of the map $\g^* \rightarrow \g$ given by contraction with $[X,\pi] \in \wedge^2\g$
is a coisotropic subalgebra of $\g$.
\end{thmno}

 We remark that the coisotropic subalgebras that arise as in the theorem
are all even dimensional, therefore they are by no means all coisotropic subalgebras. 
Using this theorem we  produce in a straightforward way families of coisotropic subalgebras  when $\g$ is one of the four classical simple complex Lie algebras or one of their split real forms.

Coisotropic subalgebras give rise to lagrangian subalgebras of the Drinfeld double $\g\oplus \g^*$ (hence also  to
 Poisson homogeneous spaces \cite{Dr226227}) 
via $\kk\mapsto \kk \oplus \kk^{\circ}$. $\cL(\g\oplus \g^*)$,
the variety of lagrangian
subalgebras of $\g\oplus \g^*$, can be endowed with a Poisson structure \cite{LuEvens1}.
It would be interesting to characterize the points of $\cL(\g\oplus \g^*)$
which correspond to the coisotropic subalgebras we constructed. Notice that $\g\oplus \g^*$ is isomorphic to   the direct sum Lie algebra $\g\oplus\g$ studied in   \cite{LuEvens2} (see Remark \ref{elu}). A further reason why coisotropic subalgebras   are interesting is that they have a counterpart in the Hopf algebra setting after quantization \cite{Nicola}.


    
Even though the above theorem is phrased entirely in terms of the Lie bialgebra $\g$, its proof involves the Poisson Lie group $G$ integrating $\g$. The paper is organized as follows. 
In Section \ref{poiac} for each $g\in G$ we consider $\h^g$, the left translation to the identity of $T_g\cO$, where $\cO$ denotes the symplectic leaf through $g$.
\emph{If} $\h^g$ is a Lie subalgebra of $\g$ then it is automatically  a coisotropic subalgebra. In Section \ref{plgs} we restrict our attention to Lie bialgebras arising from $r$-matrices and elements $g$ of the form $exp(X)$, proving the theorem stated above. Section \ref{ssla} is 
devoted to explicit examples in which $\g$ is a semi-simple Lie algebra.
In the Appendix   we present the geometric motivation that lead to considering the subspaces $\h^g$, namely pre-Poisson maps.
\\

{\bf Acknowledgments:} I learnt the simple proof of Prop. \ref{coisoalg} from  Jiang-Hua Lu. The connection to the work of Evens and Lu established in Remark \ref{elu} was suggested by the referee.
 I thank Camille Laurent and Jiang-Hua Lu for helpful conversations. I am indebted to Alberto Cattaneo and to Francesco Bonechi for remarks that improved the final version of this note. Thanks to Philippe Monnier  for a visit to Toulouse in October 2008 that helped complete this work.

\section{Coisotropic subalgebras}\label{poiac}

We recall some notions from the theory of Poisson Lie groups; we refer to the expositions \cite{Lu,KS,Ping} for more details.
 
Recall that a \emph{Poisson manifold} is a manifold $P$   endowed with a bivector field
$\Lambda\in\Gamma(\wedge^2 TP)$ satisfying $[\Lambda,\Lambda]=0$, where $[\bullet,\bullet]$ denotes the
Schouten bracket on multivector fields. 
We denote by   ${\Lambda}^{\sharp} \colon  T^*P \rightarrow TP$ the map
given by contraction with $\Lambda$.

\begin{defi}A \emph{Poisson Lie group} is a Lie group $G$ equipped with a Poisson bivector $\Lambda$ such that the multiplication map $m \colon G \times G \rightarrow G$ is a Poisson map, or equivalently such that
\begin{equation}\label{mult}
\Lambda(gh)=(L_g)_*\Lambda(h)+(R_h)_*\Lambda(g) \text{      for all } g,h \in G.
\end{equation}
\end{defi}

To every element $g$ of the Poisson Lie group $G$ we associate a \emph{subspace} of its Lie algebra $\g$ as follows:
\begin{equation}\label{h}
\h^g:=(\eta^g)^{\sharp}\;\g^*,
\end{equation}
where we use the short-hand notation 
\begin{equation}\label{eta}
\eta^g:=(L_g)_*\Lambda(g^{-1})\in \wedge^2\g.
\end{equation}
The subspace 
$\h^g$ is the left-translation to the identity of $T_{g^-1}\cO$, where $\cO$ denotes the symplectic leaf of $(G,\Lambda)$ through $g^{-1}$; in particular it is always even dimensional.
Notice that  $(\eta^g)^{\sh}\colon \g^* \to \g$ satisfies the identity
$$ (L_g)_*\circ (\Lambda(g^{-1}))^{\sh}=(\eta^g)^{\sh}\circ (L_{g^{-1}})^*.$$

\begin{defi}[{\cite[Sec. 3.1]{Lu}}]
Let $\g$ be a Lie bialgebra. A Lie subalgebra $\h$ of $\g$ is called \emph{coisotropic}\footnote{A  Lie subalgebra $\h$ is coisotropic iff the connected subgroup $H$ integrating it is a  coisotropic subgroup of $(G,\Lambda)$ (see for instance \cite{Nicola}).

Another equivalent characterization of the fact that  $\h$ is a coisotropic
Lie subalgebra is the following: $\h$ is a coisotropic submanifold of $\g$, endowed with the linear Poisson structure induced by the Lie algebra $\g^*$, and    $\h^{\circ}$ is a coisotropic submanifold of the linear Poisson manifold $\g^*$.} if its annihilator $\h^{\circ}$ is a Lie subalgebra of $\g^*$.
\end{defi}

\begin{prop}\label{coisoalg}
Let $G$ be a Poisson Lie group and $g\in G$. If $\h^g\subset \g$ is a Lie subalgebra then it is automatically a coisotropic subalgebra.
\end{prop} 
\begin{proof}
Recall that, for  every Poisson manifold $(P,\Lambda)$, there is a Lie bracket\footnote{Indeed, $T^*P$ with this bracket and the bundle map $\Lambda^{\sharp} \colon T^*P \rightarrow TP$ forms a Lie algebroid \cite{CW}.} on the space of 1-forms, inducing a Lie algebra structure on 
$(T_p\cO)^{\circ}$ for each $p\in P$ (here $\cO$ denotes the symplectic leaf through $p$). It is known that the space of left-invariant 1-forms on the Poisson Lie group $G$ is closed with respect to this bracket, and that evaluation at $e\in G$ is a Lie algebra isomorphism onto the Lie algebra $\g^*$ \cite[Sect. 2.5]{Lu}. In particular $(L_{g^{-1}})^* \colon (T_{g^{-1}}\cO)^{\circ} \rightarrow \g^*$ is a Lie algebra homomorphism, with image $(\h^g)^{\circ}$. Hence $(\h^g)^{\circ}$ is a Lie subalgebra of $\g^*$.
\end{proof}

It would be interesting to study the set $\{g\in G: \h^g \text{ is a Lie subalgebra}\}$. It is closed under inversion but is not a subgroup of $G$ (see Remark \ref{nosubgroup}). 

\begin{remark}
We are indebted to Jiang Hua Lu for pointing out the above simple proof of Prop. \ref{coisoalg}. In Appendix \ref{prepmaps} we present another proof, based on properties of the left translation $L_{g}$.
\end{remark}



\section{Poisson Lie groups arising from $r$-matrices}\label{plgs}

Let $(G,\Lambda)$ be a Poisson Lie group.
In this section we determine elements $g \in G$ for which the subspace $\h^g\subset \g$ of eq. \eqref{h}
is a Lie subalgebra, for  Prop. \ref{coisoalg} tells us that then it is a coisotropic subalgebra.

\begin{lemma}\label{liesub}
If $[\eta^g,\eta^g]=0\in \wedge^3 \g$ then $\h^g$ 
is a Lie subalgebra of $\g$.
\end{lemma}
\begin{proof}
$[\eta^g,\eta^g]=0$ iff $\overrightarrow{\eta^g}$, the 
right-invariant  bivector on $G$ whose value at the identity is $\eta^g$, is a Poisson bivector. In that case the symplectic distribution $(\overrightarrow{\eta^g})^{\sharp} \; T^*G=\overrightarrow{\h^g}$ is involutive, and this is equivalent to $\h^g$ 
being a Lie subalgebra of $\g$.
\end{proof}

\begin{defi}
Let $\g$ be a Lie algebra. An  \emph{$r$-matrix} is an element  $\pi\in \wedge^2\g$ such that $[\pi,\pi]$ is $ad$-invariant.  
\end{defi}
It is known \cite{Dr1}
that if $\pi$ is an $r$-matrix for the Lie algebra $\g$ then  
$\Lambda:={\overleftarrow{\pi}}-
{\overrightarrow{\pi}}$ makes
 $G$, any Lie group integrating  $\g$, into a Poisson Lie group.  From now on we restrict ourselves to such Poisson Lie groups. Notice that 
from definition \eqref{eta}  we get   
\begin{equation}\label{etag2}
\eta^g=\pi-Ad_g \pi.
\end{equation}

Now we are able to state the main result of this paper.

\begin{thm}\label{rmat}
Let $G$ be  a Poisson Lie group corresponding to an $r$-matrix $\pi$, $X\in \g$, $g:=exp (X)$.
Assume that 
\begin{equation}\label{condi}
[X,[X,\pi]]=\lambda[X,\pi] \text{ for some }\lambda\in \RR. \end{equation}
Then $\h^g$ is a coisotropic subalgebra of $\g$. Further 
\begin{equation}\label{hg}
\h^g={[X,\pi]}^{\sharp}\g^*.
\end{equation}

\end{thm}
\begin{proof}
Notice that
$$Ad_{exp (X)}\pi =e^{ad_X}\pi=\pi+[X,\pi]+\frac{1}{2}[X,[X,\pi]]+\frac{1}{3!}[X,[X,[X,\pi]]]+\dots
=\pi+\frac{e^{\lambda}-1}{\lambda}[X,\pi].$$
Therefore
$$\eta^g=
 \pi-Ad_g \pi=\pi-(\pi+\frac{e^{\lambda}-1}{\lambda}[X,\pi])=-\frac{e^{\lambda}-1}{\lambda}[X,\pi].$$ 

Now we use twice the fact that  $[\pi ,[X,\pi]]=\frac{1}{2}[X,[\pi,\pi]]=0$ (by the graded Jacobi identity) to show that 
$$[[X,\pi],[X,\pi]]=[X,[\pi,[X,\pi]]]- [\pi,[X,[X,\pi]]]=0-\lambda\cdot 0=0.$$
This means that $[\eta^g,\eta^g]=0$, and by Lemma
\ref{liesub} and Prop. \ref{coisoalg} $\h^g$ is a coisotropic subalgebra. 
The last part of the theorem follows since  the function 
$\frac{e^{\lambda}-1}{\lambda}$ never vanishes.
\end{proof}

\begin{remark}
If $X\in \g$ satisfies condition \eqref{condi} then 
$\Lambda={\overleftarrow{\pi}}-
{\overrightarrow{\pi}}$ and $\overrightarrow{\eta^g}$ (or $\overleftarrow{\eta^g}$) are commuting Poisson structures on $G$. This follows at once from the computations of the proof of Thm. \ref{rmat}, noticing that $\eta^g$ is a multiple of $[X,\pi]$. Here at usual $g:=exp(X)$. 
\end{remark}

We now display two very simple examples.

\begin{ep}
Let $\g=\su(2)$, so that for a suitable basis we have $[e_1,e_2]=e_3, 
[e_2,e_3]=1,[e_3,e_1]=e_2$, and take the $r$-matrix $\pi=2e_2\wedge e_3$  as in  \cite[Ex. 2.10]{Lu}. Then the only elements of $\su(2)$ that satisfy eq. \eqref{condi}
are the multiples $X$ of $e_1$, and applying \eqref{hg} we see that they all give $\h^{exp(X)}=\{0\}$.
\end{ep}

\begin{ep}\label{sl2}
Let $\g=\sll(2,\RR)$, with basis 
$$e_1=\frac{1}{2}\begin{pmatrix}
  1    & 0   \\
 0     &  -1
\end{pmatrix},\;\;\;
e_2=\frac{1}{2}\begin{pmatrix}
  0    & 1   \\
 -1     &  0
\end{pmatrix},\;\;\;
e_3=\frac{1}{2}\begin{pmatrix}
  0    & 1   \\
 1     &  0
\end{pmatrix}.$$
Then $[e_1,e_2]=e_3, [e_2,e_3]=e_1,[e_3,e_1]=-e_2$, and
$\pi=2e_2\wedge e_3$ is an $r$-matrix  \cite[Ex. 2.9]{Lu}.
The vectors $X$ of $\sll(2,\RR)$ that satisfy eq. \eqref{condi} are exactly
those of the form 
$\alpha e_1 +\beta (e_2+e_3)$ (the upper triangular matrices) and $\alpha e_1 +\beta (e_2-e_3)$ (the lower triangular matrices). 
 Applying Thm. \ref{rmat} we obtain coisotropic subalgebras $span\{e_1,e_2-e_3\}$, 
$span\{e_1,e_2+e_3\}$ and $\{0\}$.

Using \eqref{eta} one can compute directly all the elements $g\in G=SL(2\,\RR)$ for which
$[\eta^g,\eta^g]=0$: they those of the form 
 $\left(\begin{smallmatrix} 
  a    & b   \\
 0     &  a^{-1}
\end{smallmatrix}\right)$
and
$\left(\begin{smallmatrix}
  a    & 0   \\
 c     &  a^{-1}
\end{smallmatrix}\right)$. 
 By Lemma \ref{liesub} and Prop. \ref{coisoalg} these group
elements $g$ give rise to a coisotropic subalgebra of $\g$. 
The first class of elements $g$  with $b\neq 0$ all
 give rise to $span\{e_1,e_2-e_3\}$,  the   second 
 class of elements $g$  with $c\neq 0$ all give rise to
 $span\{e_1,e_2+e_3\}$, and the diagonal matrices give rise to
the trivial subalgebra $\{0\}$, i.e. we obtain exactly the same  coisotropic subalgebras as above.
\end{ep} 

\begin{remark}\label{nosubgroup}
We show that $\{g\in G: \h^g \text{ is a Lie subalgebra}\}$ is closed under the inversion map but not under multiplication. Indeed notice that $\eta^{g^{-1}}=-Ad_{g^{-1}}\eta^{g}$ 
by \eqref{mult}, so $\h^{g^{-1}}=Ad_{g^{-1}}\h^{g}$, and since $Ad_{g^{-1}}$ is a Lie algebra isomorphism the first statement follows.

To show the second statement consider $\g=\sll(2,\RR)$ as in Example \ref{sl2}. The elements
$g=\left(\begin{smallmatrix} 
  1    & 1   \\
 0     &  1
\end{smallmatrix}\right)$,
 $h=\left(\begin{smallmatrix} 
  1    & 0   \\
 -1     &  1
\end{smallmatrix}\right)$
of $G=SL(2,\RR)$ have the property that $\h^g$ and $\h^h$ are Lie subalgebras, by 
Example \ref{sl2}. However 
$\eta^{gh}=\pi-Ad_{gh}\pi=2(e_1\wedge e_2+2e_2\wedge e_3-e_1\wedge e_3)$, implying that 
$\h^{gh}$ is not a Lie subalgebra of $\g$.
\end{remark}

\section{Examples: semi-simple complex Lie algebras }\label{ssla} 

In this section we consider the standard Lie bialgebra structure on a semi-simple \emph{complex} Lie algebra, and out of its roots, using Thm. \ref{rmat} we construct families of coisotropic subalgebras. We write down explicitly\footnote{One reason for doing this is that we were not able to find any  explicit families of examples of coisotropic subalgebras  in the literature.} the resulting families for the classical simple Lie algebras $\sll(n+1,\CC),\so(2n+1,\CC),\spp(2n,\CC),\so(2n,\CC)$ and for their split real forms $\sll(n+1,\RR),
 \so(n+1,n),  
\spp(2n,\RR), \so(n,n)$. We refer to \cite[Ch. 2.6]{Ar}, to \cite{FH} and to \cite{Kn} for  background material about semi-simple complex Lie algebras and their real forms.

Let $\g$  be a semi-simple Lie algebra over $\CC$,  and fix a Cartan subalgebra $\h$. There is a decomposition $\g=\h\oplus_{\alpha\in R}\g^{\alpha}$ where
 $\g^{\alpha}$ denotes the one dimensional eigenspace for the adjoint action of $\h$ associated to the ``eigenvalue'' $\alpha\in \h^*$. The set
$R\subset \h^*$ is called   root system; make a choice $R_+$ of positive roots. 
For each $\alpha\in R_+$ choose non-zero $e_{\alpha}\in \g^{\alpha}$ and $f_{\alpha}\in \g^{-\alpha}$.

  Then an $r$-matrix is given by\begin{equation}\label{pingr}
\pi:=\sum_{{\alpha}\in R_+}\lambda_{\alpha} e_{\alpha}\wedge f_{\alpha}
\end{equation}
where $\lambda_{\alpha}:=\frac{1}{B(e_{\alpha},f_{\alpha})}$  \cite[Ex. 2.10]{Ping}. 
Notice that, since the subspaces $\g^{\alpha}$ are one dimensional and the Killing form  $B$ is $\CC$-bilinear, the above $r$-matrix depends only on the choice of Cartan subalgebra.

\begin{remark}\label{elu}
 As above let $\g$  be a semi-simple complex Lie algebra.
Evens and Lu \cite{LuEvens2}\cite[Sec. 2.1]{ELGro} consider the direct sum Lie algebra $\g\oplus\g$  endowed with the   pairing\footnote{They actually consider any non-zero multiple of the Killing form, not just $\frac{1}{2}$.}
  $\langle x_1+y_1,x_2+y_2 \rangle=\frac{1}{2}B(x_1,y_1)-\frac{1}{2}B(x_2,y_2)$ where $B$ is the Killing form of $\g$. They study the variety $\cL(\g\oplus \g)$ of lagrangian subalgebras, and endow it with interesting Poisson structures. 

Since $(\g,[\pi,\bullet])$ is a Lie bialgebra, 
$\g\oplus \g^*$ admits a Lie algebra structure known as Drinfeld double, for which the natural pairing is $ad$-invariant \cite[Sec. 2.3]{Lu}. If $\kk\subset \g$ is a coisotropic subalgebra, then $\kk \oplus \kk^{\circ}\subset \g\oplus \g^*$ is a lagrangian subalgebra.

 There is an isomorphism of Lie algebras  
\begin{equation}\label{gg}
\g\oplus \g^* \cong \g\oplus \g
\end{equation}
 preserving the pairings. As a consequence, coisotropic subalgebras of $\g$ give rise to  points of $\cL(\g\oplus \g)$, which as seen above is an interesting and well-studied variety.

Eq. \eqref{gg} follows from \cite[Prop. 1.5]{QuantumR} (see also \cite[Prop. 2.1]{Milen}). We reproduce the proof for completeness.
Recall that a Manin triple consists of a Lie algebra with an $ad$-invariant non-degenerate symmetric pairing and a decomposition into two Lagrangian subalgebras.  
There is a bijection between  Manin triples and Lie bialgebras \cite[Thm. 2.3.2]{KorSoib}. $\g\oplus \g$, together with  the diagonal $\g_{\Delta}$ and 
 \begin{equation}\label{gst}
 \{(h+v,-h+w):h\in \h, v\in \oplus_{\alpha\in R_+}\g^{\alpha}, w\in \oplus_{\alpha\in R_+}\g^{-\alpha}\},
\end{equation}
forms a Manin triple. The corresponding Lie bialgebra consists of the Lie algebra $\g$ with the derivation of $\wedge^{\bullet} \g$ obtain dualizing the Lie bracket on \eqref{gst}. A computation shows that this derivation is exactly $[\pi,\bullet]$. Hence the Drinfeld double $\g\oplus \g^*$ of the  Lie bialgebra $(\g,[\pi,\bullet])$ is isomorphic to $\g\oplus \g$ by a pairing-preserving map, showing \eqref{gg}.
\end{remark}


\begin{lemma}\label{vanish}
Let $X\in \g$ and assume that for all ${\alpha}\in R_+$
\begin{itemize}
\item[1)]$[X,[X,e_{\alpha}]]\wedge f_{\alpha}=0$
\item[2)] $[X,[X,f_{\alpha}]]\wedge e_{\alpha}=0$
\item[3)]$[X,e_{\alpha}]\wedge [X,f_{\alpha}]=0.$
\end{itemize}
Then $X$ satisfies  condition  \eqref{condi} (with $\lambda=0$).
\end{lemma}
\begin{proof}
We compute 
$$[X,\pi]=\sum_{{\alpha}\in R_+}\lambda_{\alpha} ([X,e_{\alpha}]\wedge f_{\alpha}
+  e_{\alpha}\wedge [X,f_{\alpha}]),$$ so
$$[X,[X,\pi]]=\sum_{{\alpha}\in R_+}\lambda_{\alpha} ([X,[X,e_{\alpha}]]\wedge f_{\alpha}
+2 [X,e_{\alpha}]\wedge [X,f_{\alpha}]
+  e_{\alpha}\wedge [X[X,f_{\alpha}]]),
$$ each term of which vanishes by our assumptions. 
\end{proof}

\begin{prop}\label{propline}
Let $\beta \in R_+$ satisfy this condition:
\begin{equation}\label{line}
\text{For all }\alpha \in R:\;\;\;\;(\alpha+\ZZ\beta)\cap R\text{ does not contain a string
of 3 consecutive elements.}
\end{equation}
Then $e_{\beta}$ and $f_{\beta}$ satisfy   condition  \eqref{condi}.
\end{prop}
 
\begin{proof}
We check that $X=e_{\beta}$ satisfies the assumptions of Lemma \ref{vanish}; the proof for $f_{\beta}$ is similar. Let $\alpha\in R$.

Suppose that $[e_{\beta},[e_{\beta},e_{\alpha}]]\neq 0$. Then $\alpha, \alpha+\beta$ and 
$\alpha+2\beta$
form a string of 3 consecutive elements
in $(\alpha+\ZZ\beta)\cap (R\cup \{0\})$.
Since the intersection of $R$ with any line through the origin   is
either empty or of the form $\{\alpha,-\alpha\}$ \cite[Prop. 2.20]{Ar} it follows that
$\beta=-\alpha$. So $[e_{\beta},[e_{\beta},e_{\alpha}]]$ is a multiple of $f_{\alpha}$, and assumption 1) of
 Lemma \ref{vanish} is satisfied.

Similarly, if $[e_{\beta},[e_{\beta},f_{\alpha}]]\neq 0$, then  $-\alpha, -\alpha+\beta$ and 
$-\alpha+2\beta$
form a string of 3 consecutive elements
in $(\alpha+\ZZ\beta)\cap (R\cup \{0\})$, so we must have
$\beta=\alpha$. So $[e_{\beta},[e_{\beta},f_{\alpha}]]$ is a multiple of $e_{\alpha}$, and   assumption 2) of
 Lemma \ref{vanish} is satisfied.

At most one of $\alpha+\beta$ or $\alpha-\beta$ 
 lie in $R$: if they both did then $\{\alpha-\beta, \alpha, \alpha+\beta\}$ would be a string of 
  3 consecutive elements in  $(\alpha+\ZZ\beta)\cap R$, contradicting our assumption. 
  If $\alpha-\beta \notin R$ then either $\alpha-\beta=0$, in which case $[e_{\alpha},e_{\beta}]=0$, or $[e_{\alpha},f_{\beta}]\in \g^{\alpha-\beta}=\{0\}$. A similar reasoning holds for $\alpha+\beta$, so we conclude that
  assumption 3) of
 Lemma \ref{vanish} holds.
\end{proof}

\begin{cor}
Assume the notation above and assume that
$\beta\in R_+$ satisfy  condition \eqref{line}. Let $\g_{\RR}$  denote $\g$ viewed as a real Lie algebra.  Then 
$[e_{\beta},\pi]^{\sharp}{\g_{\RR}}^*$ and $[f_{\beta},\pi]^{\sharp}{\g_{\RR}}^*$
\begin{itemize}
\item  are coisotropic subalgebras of $\g_{\RR}$
\item their complexifications
 are coisotropic subalgebras of the complex Lie bialgebra $\g$.
\end{itemize}
\end{cor}
\begin{proof}
The first statement follows from Prop. \ref{propline} and applying Thm. \ref{rmat} to $\g_{\RR}$.

 Now choose $\tilde{e}_{\alpha}\in \g^{\alpha}$ and $\tilde{f}_{\alpha}\in \g^{-\alpha}$ to be part of a Chevalley 
basis  \cite[Ch. 2.6]{Ar} of $\g$, so that $$\g_0:=
 \{h\in \h:\alpha(h)\in \RR \text  { for all }\alpha\in R_+\} \oplus_{\alpha\in R_+}span_{\RR}\{\tilde{e}_{\alpha},\tilde{f}_{\alpha}\}$$ is 
 a Lie subalgebra of $\g_{\RR}$, namely a split real form of $\g$ \cite[p. 296]{Kn}.
Since $\pi \in \wedge^2\g_0$ and $\tilde{e}_{\beta}\in \g_0$, applying Thm. \ref{rmat} to $\g_{0}$
we deduce that $[\tilde{e}_{\beta},\pi]^{\sharp}{\g_{0}}^*$ is a coisotropic subalgebra of $\g_0$. The complexification of $[\tilde{e}_{\beta},\pi]^{\sharp}{\g_{0}}^*=
[\tilde{e}_{\beta},\pi]^{\sharp}{\g_{\RR}}^*$ coincides with the complexification of 
$[e_{\beta},\pi]^{\sharp}{\g_{\RR}}^*$, hence the second statement follows.
\end{proof}

Our main references for the computation of the examples below are \cite[part III]{FH} and \cite{Var}. Two remarks about   the derivation of the examples are in order.
\begin{remark}

1) We use the fact that the Killing form $B(A_1,A_2)$   is a non-zero real multiple of $Tr(A_1A_2)$  \cite[Ex. 14.36]{FH}. 
Since the elements $e_{\alpha}$ and $f_{\alpha}$ we choose are always \emph{real} matrices,  the bivector $\pi$ is also real, and  the coisotropic subalgebras of $\g_{\RR}$   we obtain are also coisotropic subalgebras of $\g\cap Mat(n,\RR)$, which agrees with the split real form of $\g$.

2) The coisotropic subspace associated to $f_{\beta}$ will be obtained just applying the transposition map to the one associated to $e_{\beta}$. 
Indeed in all the examples below the transposition map $\bullet^T$ is an anti-homomorphism of $\g$ which switches the $e_{\alpha}$'s and the $f_{\alpha}$'s, so it maps $\pi$ to $-\pi$ and $[e_{\beta},\pi]$ to $[f_{\beta},\pi]$.
\end{remark}

\begin{ep}[$A_n$]\label{An}
Let $\g=\sll(n+1,\CC)$  with Cartan subalgebra $\h$ given by the
diagonal matrices, so that as roots we obtain
  $R=\{L_i-L_j\}_{(i\neq j)}\subset \RR^{n+1}$, where $L_1,\cdots,L_{n+1}$ denotes the standard basis of $\RR^{n+1}$.
  It is easy to  check that all roots satisfy   assumption \eqref{line}.

For a root $\alpha=L_i-L_j$ with $i<j$ we  choose $e_{\alpha}:=E_{ij}\in \g^{L_i-L_j}$ and 
$f_{\alpha}:=E_{ji}\in \g^{-L_i+L_j}$,  where 
$E_{ij}$ denotes the matrix with $1$ in the $(i,j)$-entry and zeros
elsewhere. 
We have
$\pi \sim\sum_{i<j}E_{ij}\wedge E_{ji}$, where ``$\sim$'' means ``is a non-zero real multiple of''. Fix a  root $\beta= L_{i}-L_{j}$ with $i<j$. A computation shows that 
$$[E_{ij}, \pi]\sim 
\big(\sum_{i<k\le j}+\sum_{i\le k< j}\big)E_{ik}\wedge E_{kj}=
2\sum_{i<k<j}E_{ik}\wedge E_{kj} 
-E_{ij}\wedge (H_{i}-H_{j}),$$ 
where $H_i:=E_{ii}$,
so for all $i<j$ we obtain a coisotropic subalgebra of $\g$ spanned by
\begin{center}
\fbox{$E_{ij},\;\;\;\;H_{i}-H_{j},\;\;\;\;\{E_{kj}\}_{i<k<j} \text{ and }\{E_{ik}\}_{i<k<j}$}. \end{center}

For instance, letting $n=2$ and taking $e_{\beta}=E_{13}$ leads to the
coisotropic subalgebra 
$$\left\{ \begin{pmatrix}
  a    & b &c   \\
 0     &  0 & d \\
 0     &  0 & -a  \\
\end{pmatrix}: a,b,c,d \in \RR
\right\}.
$$
 
The coisotropic subalgebra
we obtain from 
$f_{\beta}=E_{ji}$ ($i<j$)
is spanned by
\begin{center}
\fbox{$E_{ji},\;\;\;\;H_{i}-H_{j},\;\;\;\;\{E_{ki}\}_{i<k<j} \text{ and }\{E_{jk}\}_{i<k<j}$}. \end{center}

All of the above are also  coisotropic subalgebras of the split real form  $\sll(n+1,\RR)$.
\end{ep}

\begin{ep}[$B_n$]\label{Bn}
Let $\g=\so(2n+1,\CC)$,
with Cartan subalgebra given by the
diagonal matrices.
Then $R=\{\pm L_i\pm L_j\}_{(i<j)}\cup \{\pm L_i\} \subset \RR^n$.
The roots that satisfy   assumption \eqref{line} are exactly those of the
form $\pm L_i\pm L_j$ ($i< j$).

The root space of a root $L_i-L_j$ (with $i\neq j$) is spanned by $X_{ij}=E_{i,j}-E_{n+j,n+i}$.
The root space of a root $L_i+L_j$ is spanned by $Y_{ij}=E_{i,j+n}-E_{j,n+i}$,
 the one of
$-L_i-L_j$ is spanned by $Z_{ij}=E_{n+i,j}-E_{n+j,i}$. Finally,
the root space  of $L_i$ is spanned by $U_{i}=E_{i,2n+1}-E_{2n+1,n+i}$ and the
one   of $-L_i$ is spanned by $V_{i}=E_{n+i,2n+1}-E_{2n+1,i}$. As earlier, $E_{ij}$ denotes the matrix with $1$ in the $(i,j)$-entry and zeros
elsewhere. 
The $r$-matrix of eq. \eqref{pingr} satisfies $$\pi \sim\frac{1}{2}\big(\sum_{i<j}X_{ij}\wedge X_{ji}- \sum_{i<j}Y_{ij}\wedge Z_{ij}
- \sum_{i}U_{i}\wedge V_{i}\big).$$

Given a root $\beta=L_{i}-L_{j}$ (with $i<j$), a lengthy but straightforward computation shows
$$[X_{ij}, \pi]\sim -2\sum_{i<k<j}\big(X_{ik}\wedge X_{kj}\big)+X_{ij}\wedge(H_{i}-H_{j}).$$ So for all $i<j$ 
we obtain a coisotropic subalgebra  spanned by
\begin{center}
\fbox{$\{X_{ik}, X_{kj}\}_{(i<k<j)},\;\;\; X_{ij},\;\;\;H_{i}-H_{j}$}\end{center}
where $H_{i}:=E_{i,i}-E_{n+i,n+i}\in \h$.
The negative root vector
$f_{\beta}=X_{ji}$ delivers the coisotropic subalgebra  spanned by
\begin{center}
\fbox{$\{X_{ki}, X_{jk}\}_{(i<k<j)},\;\;\; X_{ji},\;\;\;H_{i}-H_{j}$}.\end{center}

If instead we pick a root $\beta=L_{i}+L_{j}$  (with $i<j$) we obtain
$$[Y_{ij}, \pi]= -2\sum_{i<k\neq j} (X_{ik}\wedge Y_{kj} )
+2 \sum_{j<k} (X_{jk}\wedge Y_{ki})
+Y_{ij}\wedge(H_{i}-H_{j})
+2U_{i}\wedge U_{j},$$ 
giving rise to a coisotropic subalgebra  spanned by
\begin{center}
\fbox{$\{X_{ik},Y_{kj}\}_{(i<k\neq j)},\;\;\; \{X_{jk},Y_{ki}\}_{(j<k)},\;\;\;
Y_{ij},\;\;\;
H_{i}-H_{j},\;\;\;
U_{i},\;\;\; U_{j}$}.\end{center}  The  root $-(L_{i}+L_{j})$  (with $i<j$) delivers the 
  coisotropic subalgebra  spanned by
\begin{center}
\fbox{$\{X_{ki},Z_{kj}\}_{(i<k\neq j)},\;\;\; \{X_{kj},Z_{ki}\}_{(j<k)},\;\;\;
Z_{ij},\;\;\;
H_{i}-H_{j},\;\;\;
V_{i},\;\;\; V_{j}$}.\end{center}

All of the above are also coisotropic subalgebras of the split real form  $\so(n+1,n)$. 
\end{ep}

\begin{ep}[$C_n$]\label{Cn}
Let $\g=\spp(2n,\CC)$. Then, choosing the diagonal
matrices as Cartan subalgebra, $R=\{\pm L_i\pm L_j\}\subset \RR^n$.
The only roots that satisfy   assumption
 \eqref{line} are those of the form $\pm 2L_i$.

For $i\neq j$ the root space of a root $L_i-L_j$ is spanned by $X_{ij}=E_{i,j}-E_{n+j,n+i}$,
as in Ex. \ref{Bn};
the root space of a root $L_i+L_j$ is spanned by $Y_{ij}=E_{i,n+j}+E_{j,n+i}$,
 the one of
$-L_i-L_j$ is spanned by $Z_{ij}=E_{n+i,j}+E_{n+j,i}$. Finally,
the root space  of $2L_i$ is spanned by $U_{i}=E_{i,n+i}$ and the
one   of $-2L_i$ is spanned by $V_{i}=E_{n+i,i}$. 
We obtain the $r$-matrix
$$\pi \sim\frac{1}{2}\sum_{i<j}X_{ij}\wedge X_{ji}+\frac{1}{2}\sum_{i<j}Y_{ij}\wedge Z_{ij}
+\sum_{i}U_{i}\wedge V_{i}.$$

Let us consider the root $2L_{i}$.
A computation shows
$$[U_{i} ,\pi]\sim\sum_{i<k}(Y_{ik}\wedge X_{ik})+U_{i}\wedge
H_{i},$$
where $H_{i}:=E_{ii}-E_{n+i,n+i},$
so as coisotropic subspace we obtain the span of
\begin{center}
\fbox{$\{Y_{ik},X_{ik} \}_{i<k},\;\;\;
U_{i},\;\;\;H_{i}$}.\end{center}

For instance, when $n=2$, taking $e_{\beta}=U_2=E_{24}$ and
$e_{\beta}=U_1=E_{13}$
we obtain the coisotropic subalgebras of $\spp(4,\CC)$
$$\left\{ \begin{pmatrix}
  0    & 0 &0 & 0   \\
 0     &  a & 0 &b \\
 0     &  0 & 0  &0 \\
  0     &  0 & 0 &-a \\
\end{pmatrix}: a,b \in \RR
\right\} \text{     and      }
\left\{ \begin{pmatrix}
  a    & c &b & d   \\
 0     &  0 & d &0 \\
 0     &  0 & -a  &0 \\
  0     &  0 & -c &0 \\
\end{pmatrix}: a,b,c,d \in \RR
\right\}.$$

For the root $-2L_{i}$, whose root space is spanned by $V_{i}$,
 as coisotropic subspace we obtain the span of
\begin{center}
\fbox{$\{Z_{ik},X_{ki} \}_{i<k},\;\;\;
V_{i},\;\;\;H_{i}$}.\end{center}
All of the above are also coisotropic subalgebras of the split real form  $\spp(2n,\RR)$. 
\end{ep}

\begin{ep}[$D_n$] \label{Dn}
Let $\g=\so(2n,\CC)$. Then $R=\{\pm L_i\pm L_j\}_{\{i<j\}} 
\subset \RR^n$, and the same computation as in Ex. \ref{Bn} shows that all
  roots   satisfy   assumption
 \eqref{line}. The root spaces of $L_i-L_j,L_i+L_j$ and $-L_i-L_j$ are spanned
 by elements $X_{ij}, Y_{ij}$ and $Z_{ij}$ defined by the same formulae as in Ex. \ref{Bn}, and
the $r$-matrix of eq. \eqref{pingr} satisfies
$$\pi \sim\frac{1}{2}\big(\sum_{i<j}X_{ij}\wedge X_{ji}- \sum_{i<j}Y_{ij}\wedge Z_{ji})$$
(it consists of the first two summands of the $r$-matrix for the
 $B_n$ case).

The same computations as  in Ex. \ref{Bn} show that (with $i<j$)
from the root $L_{i}-L_{j}$ 
we obtain the coisotropic subalgebras  spanned by
 
\begin{center}
\fbox{$\{X_{ik}, X_{kj}\}_{(i<k<j)},\;\;\; X_{ij},\;\;\;H_{i}-H_{j}$}\end{center} and
\begin{center}
\fbox{$\{X_{ki}, X_{jk}\}_{(i<k<j)},\;\;\; X_{ji},\;\;\;H_{i}-H_{j}$},\end{center}
whereas from the root $L_{i}+L_{j}$   we obtain
the coisotropic subalgebras   spanned by
\begin{center}
\fbox{$\{X_{ik},Y_{kj}\}_{(i<k\neq j)},\;\;\; \{X_{jk},Y_{ki}\}_{(j<k)},\;\;\;
Y_{ij},\;\;\;
H_{i}-H_{j}$}\end{center} and
\begin{center}
\fbox{$\{X_{ki},Z_{kj}\}_{(i<k\neq j)},\;\;\; \{X_{kj},Z_{ki}\}_{(j<k)},\;\;\;
Z_{ij},\;\;\;
H_{i}-H_{j}$}.\end{center}
(Here  $H_{i}:=E_{i,i}-E_{n+i,n+i}$).
All of the above are also  coisotropic subalgebras of the real form  $\so(n,n)$. 
\end{ep}

\begin{remark}\label{notid}
In Example \ref{An}, taking $n=2$ and $g=exp (E_{13})$, we showed that  $\h^g=span_{\RR}
\{E_{12},E_{13},E_{23},H_1-H_3\}$ is a coisotropic subalgebra of $\sll(3,\RR)$. 
In particular its annihilator $(\h^g)^{\circ}$ is a Lie subalgebra, but it is \emph{not} a Lie ideal. Indeed, taking the basis of $\sll(3,\RR)$ given by $\{E_{ij}\}_{(i\neq j)}$, $H_1-H_2$, $H_1-H_3$ and considering the dual basis, we have $(H_1-H_2)^* \in (\h^g)^{\circ}$ but
$\langle [(E_{12})^*, (H_1-H_2)^* ], E_{12} \rangle \neq 0$.
\end{remark}

\appendix
\section{Pre-Poisson maps}\label{prepmaps}

In this appendix we  
generalize the notion of Poisson map between Poisson manifolds. A  natural example is the left  translation $L_g$ on a Poisson Lie group $G$ (Lemma \ref{decifit}), which gives rise naturally to the subspace $\h^g\subset T_eG$ considered in Section \ref{poiac}, providing an alternative proof of Prop. \ref{coisoalg}.

Recall that a  submanifold $C$ of a Poisson manifold $P$ is called \emph{coisotropic} if 
${\Lambda}^{\sharp}  N^*C\subset TC$, where $N^*C$ (the conormal bundle of
$C$) is defined as the annihilator of $TC$. 
Here we need a generalization of the notion of coisotropic submanifold:

\begin{defi}
A submanifold $C$ of a Poisson manifold $(P,\Lambda)$ is called
\emph{pre-Poisson} \cite{CZ} if the rank of $TC+{\Lambda}^{\sharp}  N^*C$ is constant
along $C$, or equivalently if  $ pr_{NC}\circ{\Lambda}^{\sharp} \colon N^*C\rightarrow TP|_C \rightarrow NC:=TP|_C/TC$ has constant rank.

A map $\phi \colon (P_1,\Lambda_1)\rightarrow (P_2,\Lambda_2)$ between Poisson manifolds is a  \emph{pre-Poisson map} if $graph(\phi)$ is a pre-Poisson submanifold of the product $P_1 \times \bar{P_2}$, where  $\bar{P_2}$ denotes the Poisson manifold $(P_2,-\Lambda_2)$.
\end{defi} 

A map between Poisson manifolds is a Poisson map iff its graph is coisotropic, hence we see that pre-Poisson maps generalize the notion of Poisson map. 
We make more explicit what it means to be a pre-Poisson map.
\begin{lemma}\label{const}
A map $\phi \colon (P_1,\Lambda_1)\rightarrow (P_2,\Lambda_2)$ is pre-Poisson iff for all $x\in P_1$ the rank of 
\begin{equation*} \label{nondi} E(x)=\{{(\Lambda_2-\phi_*\Lambda_1)}^{\sharp}\xi:\xi \in T^*_{\phi(x)}P_2\}
\subset T_{\phi(x)}P_2
\end{equation*} is constant. Here $\phi_* \colon T_xP_1 \rightarrow T_{\phi(x)}P_2$.\end{lemma}

\begin{proof}
Let $\Gamma:=graph(\phi)\subset P_1\times \bar{P}_2$ and $x\in P_1$. We have
\begin{eqnarray*}
T_{(x,\phi(x))}\Gamma+{(\Lambda_1-\Lambda_2)}^{\sharp}N_{(x,\phi(x))}^*\Gamma&=&\{(X,\phi_*X):X\in T_xP_1\}+\{({\Lambda}^{\sharp}_1\phi^*\xi, {\Lambda}^{\sharp}_2\xi):\xi\in T^*_{\phi(x)}P_2\}\\
&=&\{(X,\phi_*X):X\in T_xP_1\}+\{(0,{\Lambda}^{\sharp}_2\xi-\phi_*({\Lambda}^{\sharp}_1\phi^*\xi)):\xi\in T^*_{\phi(x)}P_2\}\\
&=&\{(X,\phi_*X):X\in T_xP_1\}+\{0\}\times E(x).
\end{eqnarray*}
A complement of this subspace in $T_{(x,\phi(x))}(P_1\times P_2)$ is $(0,R(x))$, where
$R(x)$ is a complement to $E(x)$ in $T_{\phi(x)}P_2$. Hence $\Gamma$ is a pre-Poisson submanifold iff $R(x)$, or equivalently $E(x)$, has constant rank as $x$ varies over all points of $P_1$.
\end{proof}

\begin{remark}\label{sur}
1) The composition of pre-Poisson maps is \emph{not} pre-Poisson. Let $P_1=(\RR^2, \frac{\partial}{\partial x}\wedge 
 \frac{\partial}{\partial y})$, $P_2=(\RR^2, 0)$ and 
 $P_3=(\RR^2, (1+x^2+y^2) \frac{\partial}{\partial x}\wedge 
 \frac{\partial}{\partial y})$. The identity maps $id \colon P_1\rightarrow P_2$ and 
 $id \colon P_2\rightarrow P_3$ are pre-Poisson maps (this is seen easily using Lemma \ref{const}), however the composition is not.

 2) Let $P_1,P_2$ be Poisson manifolds and 
$\phi \colon P_1 \rightarrow P_2$ be a submersive \emph{Poisson} map.  
 If $C\subset P_2$ is a pre-Poisson submanifold (for example a point), then $f^{-1}(C)$ is a pre-Poisson submanifold of $P_1$ \cite{CZbis}. 
When $\phi$ is just a submersive \emph{pre-Poisson} map this  statement is not longer true:
the projection  $\phi:(\RR^3, -z^2\frac{\partial}{\partial x}\wedge 
 \frac{\partial}{\partial y})\rightarrow (\RR^2, \frac{\partial}{\partial x}\wedge 
 \frac{\partial}{\partial y})$ onto the first two components  is a pre-Poisson map, but
 $\phi^{-1}(0)=\{(0,0,z):z\in \RR \}$ is not a pre-Poisson submanifold.
 \end{remark}

From now on we consider only the case when 
 the map $\phi$ of Lemma \ref{const} is a \emph{diffeomorphism}. Then  $D_y:=E(\phi^{-1}(y))$ defines a singular distribution on $P_2$ which  measures how $\phi$ fails to be a Poisson map.
\begin{defi}\label{defidef}
Given a diffeomorphism $\phi \colon (P_1,\Lambda_1)\rightarrow (P_2,\Lambda_2)$ between Poisson manifolds, the \emph{deficit distribution} associated to $\phi$ is
the singular distribution on $P_2$ given by
\begin{equation*}\label{D}
D=\{{(\Lambda_2-\phi_*\Lambda_1)}^{\sharp}\xi:\xi \in T^*P_2\}.
\end{equation*}
\end{defi}

The deficit distribution 
$D$ singles out an interesting subalgebra of $C^{\infty}(P_2)$:
\begin{lemma}\label{inv}
Let  $\phi \colon (P_1,\Lambda_1)\rightarrow (P_2,\Lambda_2)$ be a diffeomorphism. Then the set of $D$-invariant functions $\{f: d_yf|_{D_y}=0 \text{ for all }y\in P_2\}$ coincides with 
\begin{equation} \label{subal} \big\{f: \phi^*\{f,g\}=\{\phi^*f,\phi^*g\}
 \text{ for all } g\in C^{\infty}(P_2)\big\},
\end{equation}
and is a Poisson subalgebra of 
 $C^{\infty}(P_2)$.
\end{lemma}
\begin{proof}
Expressing $D$ in terms of hamiltonian vector fields we have $D=
\{X^{P_2}_g-\phi_*(X^{P_1}_{\phi^* g}): g\in C^{\infty}(P_2) \}$.
The claimed equality follows from
$$d_yf( X^{P_2}_g-\phi_*(X^{P_1}_{\phi^* g})) =\{f,g\}_y-d_{\phi^{-1}(y)}(\phi^*f)X^{P_1}_{\phi^* g}=
(\phi^*\{f,g\}-
\{\phi^*f,\phi^*g\})_{{\phi^{-1}(y)}}$$
for all $y\in P_2$.

To show that   \eqref{subal} is a Poisson subalgebra we compute
for $D$-invariant functions $f$ and $\tilde{f}$ on $P_2$ and for $g\in C^{\infty}(P_2)$ that
$$\phi^*\{\{f,g\},\tilde{f}\}=\{\phi^*\{f,g\},\phi^*\tilde{f}\}=\{\{\phi^*f,\phi^*g\},\phi^*\tilde{f}\}.$$
Hence using twice the Jacobi identity we obtain
  \begin{eqnarray*} 
 \phi^*\{\{f,\tilde{f}\},g\}
&=&\phi^*\{\{f,g\},\tilde{f}\}+\phi^*\{f,\{\tilde{f},g\}\}\\
&=&\{\{\phi^*f,\phi^*g\},\phi^*\tilde{f}\}+\{\phi^*f,\{\phi^*\tilde{f},\phi^*g\}\}=
\{\{\phi^*f,\phi^*\tilde{f}\},\phi^*g\}=\{\phi^*\{f,\tilde{f}\},\phi^*g\}.
\end{eqnarray*}
\end{proof}

Summarizing the above results   we have
\begin{cor}\label{quot}
A diffeomorphism $\phi \colon (P_1,\Lambda_1)\rightarrow (P_2,\Lambda_2)$ is a pre-Poisson map
iff  $\Lambda_2-\phi_*\Lambda_1$ is a constant rank bivector on $P_2$, i.e. iff
$D$ is a smooth constant rank distribution on $P_2$. If $D$ is integrable and the leaf space $P_2/D$ is smooth, then $P_2/D$  has a Poisson structure induced by the projection map $\pi \colon P_2 \rightarrow P_2/D$. In this case the
composition $\pi\circ\phi \colon  P_1 \rightarrow P_2/D$ is a Poisson map.
\end{cor}
\begin{proof}
$\phi$ is a pre-Poisson map by Lemma \ref{const}.
By the second part of Lemma \ref{inv}
the $D$-invariant functions on $P_2$ form a Poisson subalgebra of 
 $C^{\infty}(P_2)$, so $P_2/D$ has an induced Poisson structure.
 By the first part of Lemma \ref{inv} in particular
 $\phi^*\{f,\tilde{f}\}=\{\phi^*f,\phi^*\tilde{f}\}$ for all $D$-invariant functions    $f,\tilde{f}$   on $P_2$, so 
  $\pi\circ\phi$ is a Poisson map.\end{proof}

Now let $G$ be a Poisson Lie group and $g\in G$.
The subspace $\h^g$ defined in Section \ref{poiac} generates the deficit distribution 
of the left translation $L_g \colon G \rightarrow G$.
\begin{lemma}\label{decifit}
a) $L_g \colon G \rightarrow G$ is a pre-Poisson map.\\ 
b) Its deficit distribution is $\overrightarrow{\h^g}$, the right-invariant distribution obtained translating $\h^g\subset T_eG$.
\end{lemma}
 \begin{proof}
a) By Cor. \ref{quot} we have to show that $\Lambda- (L_g)_*\Lambda$ is a constant rank bivector on $G$. This bivector field at the point $k\in G$ is
 \begin{equation}\label{difference}
 \Lambda(k)- (L_g)_*[\Lambda(g^{-1}k)]=-(L_g)_*(R_k)_*\Lambda(g^{-1})=
 -(R_k)_*\eta^g,\end{equation}
where we have used \eqref{mult} applied to $\Lambda(g^{-1}k)$ in the first equality. In other words
$\Lambda- (L_g)_*\Lambda=-\overrightarrow{\eta^g}$, which  obviously has constant rank.

b) Using a) we see that the deficit distribution is   $[\Lambda- (L_g)_*\Lambda]^{\sharp}T^*G=[\overrightarrow{\eta^g}]^{\sharp}T^*G=\overrightarrow{\h^g}$.
\end{proof}
 

The  observations above allow for an alternative, perhaps more geometric, proof of Prop. \ref{coisoalg}.
\begin{proof}[Alternative proof of Prop. \ref{coisoalg}]
For any $f_1,f_2\in C^{\infty}(G)$ and $X\in \g$ we have \cite[Ch. 2.3]{Lu}
\begin{equation}
\label{double}
 \langle [d_e f_1,d_e  f_2],X \rangle= {X}\{f_1,f_2\}.
\end{equation}
Any element of $(\h^{\g})^{\circ}$ can be realized as $d_e f$ where $f$ is a function on $G$ which is invariant along the integrable distribution 
obtained right-translating $\h^g$. This distribution coincides with the deficit distribution of $L_g \colon G \rightarrow G$ by Lemma \ref{decifit} b).
Hence, if $f_1$ and $f_2$ are invariant functions, by Lemma \ref{inv} $\{f_1,f_2\}$ is also invariant. Therefore the right hand side of \eqref{double} vanishes for all $X\in \h^g$, from which we deduce that $ [d_ef_1,d_ef_2]\in (\h^g)^{\circ}.$
\end{proof}
  

We conclude   with two remarks on Poisson actions.
\begin{remark} The considerations of Lemma \ref{decifit} can be extended to locally free  left \emph{Poisson actions} (i.e., actions for which $\sigma: G \times P \rightarrow P$ is   a Poisson map, where  $G \times P$ is equipped with the product Poisson structure).
In this case we obtain:\\
a) for all $g\in G$, $\sigma_g \colon P \rightarrow P$ is a pre-Poisson map.\\
b) the deficit distribution of $\sigma_g$ is generated by the infinitesimal  action of $\h^{g}\subset \g$.\\
 If $\h^g$ is a Lie subalgebra of $\g$ and $P/H^g$ is a smooth manifold, where
$H^g$ the connected  subgroup of $G$ integrating $\h^g$, 
then  $P/H^g$ has a Poisson structure for which the projection map $\pi \colon P\rightarrow P/H^g$ is Poisson. This is a well-known fact, see   \cite[Thm. 6]{STS} or  \cite[Prop. 3.4]{Lu}. Cor. \ref{quot} in addition tells us
that $\pi\circ \sigma_g \colon  P \rightarrow P/H^g$ is also a Poisson map.
\end{remark}

\begin{remark}
Consider the action by left multiplication  $G$ on itself, and let $g\in G$ so that $\h^g$   is a Lie subalgebra of $\g$. Then $H^g \backslash G$ (if smooth), together with the action of $G$ by right multiplication,  is a right Poisson homogeneous space (i.e., $(H^g \backslash G) \times G \rightarrow H^g \backslash G$ is a transitive right action and a Poisson map).
Further both the projection $\pi$ and 
$\pi\circ L_g \colon  G \rightarrow H^g\backslash G$ are Poisson maps which are equivariant for the $G$-actions by right multiplication. \end{remark}

\bibliographystyle{habbrv}
\bibliography{bibPLG}

\def\cprime{$'$} \def\cprime{$'$} \def\cprime{$'$} \def\cprime{$'$}
\begin{thebibliography}{10}

\bibitem{Ar}
A.~Arvanitoyeorgos.
\newblock {\em An introduction to {L}ie groups and the geometry of homogeneous
  spaces}, volume~22 of {\em Student Mathematical Library}.
\newblock American Mathematical Society, Providence, RI, 2003.
\newblock Translated from the 1999 Greek original and revised by the author.

\bibitem{CW}
A.~Cannas~da Silva and A.~Weinstein.
\newblock {\em Geometric models for noncommutative algebras}, volume~10 of {\em
  Berkeley Mathematics Lecture Notes}.
\newblock American Mathematical Society, Providence, RI, 1999.

\bibitem{CZbis}
A.~S. Cattaneo and M.~Zambon.
\newblock Pre-poisson submanifolds.
\newblock In {\em Travaux math\'ematiques.}, Trav. Math., XVII, pages 61--74.
  Univ. Luxemb., Luxembourg, 2007.

\bibitem{CZ}
A.~S. Cattaneo and M.~Zambon.
\newblock Coisotropic embeddings in {P}oisson manifolds.
\newblock {\em Trans. Amer. Math. Soc.}, 361(7):3721--3746, 2009.

\bibitem{Nicola}
N.~Ciccoli.
\newblock Quantization of co-isotropic subgroups.
\newblock {\em Lett. Math. Phys.}, 42(2):123--138, 1997.

\bibitem{Dr6871}
V.~G. Drinfel{\cprime}d.
\newblock Hamiltonian structures on {L}ie groups, {L}ie bialgebras and the
  geometric meaning of classical {Y}ang-{B}axter equations.
\newblock {\em Dokl. Akad. Nauk SSSR}, 268(2):285--287, 1983.

\bibitem{Dr1}
V.~G. Drinfel{\cprime}d.
\newblock Quantum groups.
\newblock In {\em Proceedings of the International Congress of Mathematicians,
  Vol. 1, 2 (Berkeley, Calif., 1986)}, pages 798--820, Providence, RI, 1987.
  Amer. Math. Soc.

\bibitem{Dr226227}
V.~G. Drinfel{\cprime}d.
\newblock On {P}oisson homogeneous spaces of {P}oisson-{L}ie groups.
\newblock {\em Teoret. Mat. Fiz.}, 95(2):226--227, 1993.

\bibitem{LuEvens1}
S.~Evens and J.-H. Lu.
\newblock On the variety of {L}agrangian subalgebras. {I}.
\newblock {\em Ann. Sci. \'Ecole Norm. Sup. (4)}, 34(5):631--668, 2001.

\bibitem{LuEvens2}
S.~Evens and J.-H. Lu.
\newblock On the variety of {L}agrangian subalgebras. {II}.
\newblock {\em Ann. Sci. \'Ecole Norm. Sup. (4)}, 39(2):347--379, 2006.

\bibitem{ELGro}
S.~Evens and J.-H. Lu.
\newblock Poisson geometry of the {G}rothendieck resolution of a complex
  semisimple group.
\newblock {\em Mosc. Math. J.}, 7(4):613--642, 766, 2007.

\bibitem{FH}
W.~Fulton and J.~Harris.
\newblock {\em Representation theory}, volume 129 of {\em Graduate Texts in
  Mathematics}.
\newblock Springer-Verlag, New York, 1991.
\newblock A first course, Readings in Mathematics.

\bibitem{Kn}
A.~W. Knapp.
\newblock {\em Lie groups beyond an introduction}, volume 140 of {\em Progress
  in Mathematics}.
\newblock Birkh\"auser Boston Inc., Boston, MA, second edition, 2002.

\bibitem{KorSoib}
L.~I. Korogodski and Y.~S. Soibelman.
\newblock {\em Algebras of functions on quantum groups. {P}art {I}}, volume~56
  of {\em Mathematical Surveys and Monographs}.
\newblock American Mathematical Society, Providence, RI, 1998.

\bibitem{KS}
Y.~Kosmann-Schwarzbach.
\newblock Lie bialgebras, {P}oisson {L}ie groups and dressing transformations.
\newblock In {\em Integrability of nonlinear systems (Pondicherry, 1996)},
  volume 495 of {\em Lecture Notes in Phys.}, pages 104--170. Springer, Berlin,
  1997.

\bibitem{Ping}
C.~Laurent-Gengoux, M.~Stienon, and P.~Xu.
\newblock Lectures on poisson groupoids, 2007, arXiv.org:0707.2405.

\bibitem{Lu}
J.-H. Lu.
\newblock Multiplicative and affine poisson structures on lie groups, 1990,
  Ph.D. Thesis, U.C. Berkeley, available at
  \texttt{http://hkumath.hku.hk/~jhlu/publications.html}.

\bibitem{QuantumR}
N.~Y. Reshetikhin and M.~A. Semenov-Tian-Shansky.
\newblock Quantum {$R$}-matrices and factorization problems.
\newblock {\em J. Geom. Phys.}, 5(4):533--550 (1989), 1988.

\bibitem{STS}
M.~A. Semenov-Tian-Shansky.
\newblock Dressing transformations and {P}oisson group actions.
\newblock {\em Publ. Res. Inst. Math. Sci.}, 21(6):1237--1260, 1985.

\bibitem{Var}
V.~S. Varadarajan.
\newblock {\em Lie groups, {L}ie algebras, and their representations}, volume
  102 of {\em Graduate Texts in Mathematics}.
\newblock Springer-Verlag, New York, 1984.
\newblock Reprint of the 1974 edition.

\bibitem{Milen}
M.~Yakimov.
\newblock Symplectic leaves of complex reductive {P}oisson-{L}ie groups.
\newblock {\em Duke Math. J.}, 112(3):453--509, 2002.

\end{thebibliography}
\end{document}